\documentclass[11pt,reqno]{amsart}
\usepackage[a4paper,margin=28mm]{geometry}
\usepackage{amsmath,amssymb,amsthm,mathrsfs,enumitem,microtype}
\usepackage[colorlinks=true,linkcolor=blue,citecolor=blue,urlcolor=blue]{hyperref}
\newtheorem{theorem}{Theorem}[section]
\newtheorem{proposition}[theorem]{Proposition}
\newtheorem{lemma}[theorem]{Lemma}
\newtheorem{corollary}[theorem]{Corollary}
\theoremstyle{definition}

\newtheorem{example}[theorem]{Example}

\theoremstyle{remark}
\newtheorem{remark}[theorem]{Remark}
\newcommand{\Mod}{\text{-}\mathrm{Mod}}
\newcommand{\Hom}{\operatorname{Hom}}
\newcommand{\Ext}{\operatorname{Ext}}
\newcommand{\Tor}{\operatorname{Tor}}
\newcommand{\Spec}{\operatorname{Spec}}
\newcommand{\Ann}{\operatorname{Ann}}
\newcommand{\A}{\mathcal A_I}

\title[Relative derived torsion in reduced module categories]{Relative Derived Torsion in Categories of Reduced Modules:\\ Comparison and Non-Noetherian Examples}
\author{Patrick Masaba}
\author{David Ssevviiri}
\address{Department of Mathematics, Makerere University, P.O. Box 7062, Kampala, Uganda}
\email{masabapatrick30@gmail.com}
\email{david.ssevviiri@mak.ac.ug}
\makeatletter
\renewenvironment{abstract}{%
  \global\setbox\abstractbox=\vtop\bgroup
  \normalfont\small
  \begin{quotation}
    \begin{center}
      \bfseries ABSTRACT
    \end{center}
    \setlength{\parindent}{0pt}
    \noindent\ignorespaces
}{%
  \end{quotation}
  \egroup
}
\makeatother
\begin{document}
\begin{abstract}
Let $I$ be an ideal of a commutative ring $R$. We use relative derived torsion in exact full abelian subcategories of $I$-reduced modules to study ambient derived ideal-power torsion. Exact inclusion forces $IM=I^2M$ whenever $I^2$ is finitely generated. Over a Noetherian ring, every finite object consequently decomposes into its torsion part and $IM$, and its positive-degree ambient cohomology vanishes. For $R=S[t]$ and $I=JS[t]+(t)$ with $J^2=J$, we compute every stage of the canonical Ext system and its transition maps, obtaining a comparison with $\Ext_S^q(S/J,-)$. Valuation-domain and product-of-fields examples yield nonzero cyclic first cohomology. A further construction has cyclic coefficients and cyclic ambient cohomology in every positive even degree, while its positive-degree relative cohomology vanishes. Applications determine kernels of canonical Ext maps for these coefficient modules and identify the obstruction to ring-valued specialization in the questions of Eisenbud, Musta\c{t}\u{a}, and Stillman.
\end{abstract}
\maketitle
\section{Introduction}
\subsection{Ambient local cohomology and relative resolutions}
For an ideal $I$ of a commutative unital ring $R$, and an $R$-module $M$, put
\[
\Gamma_I(M):=\bigcup_{n\geq1}(0:_M I^n),\text{ and }H_I^q(M):=R^q\Gamma_I(M).
\]
Throughout, $H_I^q$ means derived ideal-power torsion in $R\Mod$. Over Noetherian rings, these are the usual local cohomology functors. Over arbitrary rings, we do not identify them with geometric cohomology with supports without further hypotheses. The formula
\begin{equation}\label{eq:limit}
H_I^q(M)\cong\varinjlim_n\Ext_R^q(R/I^n,M)
\end{equation}
holds over any commutative ring: apply the filtered-colimit description of $\Gamma_I$ to an injective resolution and use exactness of filtered colimits. Its transition maps are essential data; a nonzero Ext class can disappear before reaching local cohomology.

We study this ambient functor through injective resolutions in a full abelian subcategory $\A\subseteq R\Mod$ with exact inclusion and enough injectives. Its objects are $I$-reduced, that is, $(0:_M I)=(0:_M I^2)$, where \((0:_M I):=\{\,m\in M \mid Im=0\,\}.\) The relative functors are obtained by deriving $\Gamma_I|_{\A}$ in $\A$, with values in $R\Mod$. The purpose is to determine when the simpler relative computation gives the ambient module $H_I^q(M)$, and to describe the obstructions when it does not.

\subsection{Which finiteness properties are being established?}
Finiteness in local cohomology has several different meanings. Over a Noetherian ring, finite generation of a local cohomology module $H_I^i(M)$ implies finiteness of its associated-prime set; the converse need not hold. Finiteness of Bass numbers is another invariant. In this classical Noetherian setting, Hartshorne's cofiniteness condition asks that a module $N$ have support in $V(I)$ and that every $\Ext_R^j(R/I,N)$ be finitely generated; it does not require $N$ itself to be finitely generated. These distinctions are fundamental in the classical theory \cite{BrodmannSharp2013,GrothendieckHartshorne1967,Hartshorne1970,Huneke1992}.

Huneke and Sharp established finite Bass numbers and associated-prime sets for local cohomology of regular local rings in positive characteristic \cite{HunekeSharp1993}. Lyubeznik obtained characteristic-zero finiteness results by $D$-module methods \cite{Lyubeznik1993}. In the cofiniteness direction, Delfino and Marley proved that local cohomology of finite modules over a Noetherian local ring is $I$-cofinite when $\dim(R/I)=1$, \cite{DelfinoMarley1997}. These results control modules which can still fail to be finitely generated over the original ring. For instance, $H_{(x)}^1(k[x])\cong k[x,x^{-1}]/k[x]$ is not finite over $k[x]$.

In this paper, we establish finite-generation results for ambient derived torsion itself. Their precise scope is as follows: For finite coefficient modules in the Noetherian categorical setting, every $H_I^q(M)$ is finite, and all positive-degree groups vanish. Outside the Noetherian setting, we construct nonzero cyclic modules $H_I^1(M)\cong R/I$ for infinitely generated coefficient ideals. We also give a cyclic coefficient module $k$ for which $H_I^{2r}(k)\cong k$ for all $r\geq1$ and the positive odd-degree groups vanish. Thus, actual finite generation can coexist with positive-degree nonvanishing in these non-Noetherian examples. Their supporting ideals are infinitely generated. The Noetherian conclusion is thus a vanishing theorem under categorical restrictions; the non-Noetherian examples provide nonzero finite-generation results.

\subsection{Relation with $I$-reduced and $I$-coreduced modules}
Ssevviiri proved the natural description $\Gamma_I(M)\cong\Hom_R(R/I,M)$ for $I$-reduced modules and developed its connections with completion and module-theoretic duality \cite{Ssevviiri2024}. The subsequent work \cite{SsevviiriII} studies radicality and torsion-theoretic correspondences on suitable subcategories. We use these degree-zero identities as input. The distinction here is the explicit choice of an abelian category in which resolutions are taken and the verification that its internal injectives compute ambient derived torsion. A degree-zero identity on reduced objects alone does not identify their ambient derived functors.

There is also a categorical restriction beyond reducedness of an individual module. Proposition~\ref{prop:coreduced} shows that exact inclusion forces $IM=I^2M$ i.e., ($M$ to be $I$-coreduced $R$-module)  on all objects whenever $I^2$ is finitely generated. This is obtained by realizing $M/I^2M$ as an internal cokernel; arbitrary quotient closure is not needed. It leads to the finite-coefficient vanishing theorem and prevents such a subcategory from containing $R$ for a non-idempotent ideal in a Noetherian ring. The distinction is visible already for $R=k[x]$: $R$ itself is $(x)$-reduced, but has nonzero $H_{(x)}^1(R)$ and cannot belong to an exact full abelian subcategory consisting entirely of $(x)$-reduced modules.

\subsection{Standard homological tools and the comparison problem}
Proposition~\ref{Comp:Thm} is an application of the acyclic-resolution theorem: relative and ambient derived torsion agree on all objects exactly when internal injectives are ambient $\Gamma_I$-acyclic. Likewise, Proposition~\ref{MainThm} uses resolutions acyclic for both $\Hom_R(R/I,-)$ and $\Gamma_I$. These standard derived-functor arguments \cite{Weibel1994} make the transfer hypotheses explicit.

Categorical local cohomology itself also predates this work. Savoji and Sazeedeh define section functors using open subsets of the atom spectrum of a locally Noetherian Grothendieck category and study their derived functors \cite{SavojiSazeedeh2019}. Here no atom-spectrum construction or locally Noetherian Grothendieck hypothesis is imposed on $\A$: the functor is the restriction of a specified ambient ideal-power torsion functor, and the question is whether its internal derivation computes the ambient one. The ambient comparison, rather than the existence of categorical derived section functors, is the issue addressed here.

For $S=R/K$, with idempotent image $J=(I+K)/K$, the subcategory of modules killed by $K$ is $S\Mod$, and the relative computation is $\Ext_S^q(S/J,M)$. Corollary~\ref{cor:quotient-comparison} gives a sufficient ambient comparison condition: the systems $\{\Tor_q^R(S,R/I^n)\}_n$ are pro-zero for $q>0$. The proof uses the usual change-of-rings adjunction into an $S$-injective module. This specializes change-of-rings theory to a usable ambient comparison condition for the quotient categories.

\subsection{Non-Noetherian torsion and existing comparison theories}
Schenzel studied the relation between local cohomology, \v{C}ech complexes, and weak proregularity \cite{Schenzel2003}. Porta, Shaul, and Yekutieli developed derived completion and torsion and the MGM equivalence under weak proregularity \cite{PortaShaulYekutieli2014}. Vyas and Yekutieli related weak proregularity to weak stability of torsion classes \cite{VyasYekutieli2018}; Yekutieli further distinguished idealistic and sequential derived constructions \cite{Yekutieli2021}. Our relative construction concerns a different comparison: the same ideal-power torsion functor derived in two different abelian categories. In particular, the infinitely generated supporting ideals in our examples are outside the usual finitely generated weak-proregularity hypotheses.

Qu\'y and Rohrer study preservation of injectivity by torsion functors, called ITI, and give non-Noetherian comparison and nonvanishing results \cite{QuyRohrer2017}. Their Proposition~3.7 already gives higher local cohomology for infinitely generated ideals in non-Noetherian absolutely flat rings; their Proposition~2.2 supplies valuation-type and zero-dimensional examples of failure of ITI. The additional information obtained here is the explicit cyclic value of the cohomology in the examples, the coefficient-module finiteness, and the full polynomial Ext transition system. ITI concerns torsion of ambient injectives, whereas our comparison tests ambient acyclicity of internal injectives.

In particular, idempotence of $J$ only gives $\Gamma_J=\Hom_S(S/J,-)$. This functor is exact precisely when $S/J$ is projective. The finitely generated idempotent case splits, but an infinitely generated idempotent ideal need not split. This standard distinction is central to the examples and prevents treating idempotence alone as an ambient vanishing criterion.

\subsection{Polynomial lifting and canonical Ext maps}
Our main explicit comparison calculation is Theorem~\ref{thm:lift}. For
\[
R=S[t],\quad I=JS[t]+(t),\quad J^2=J,\quad tM=0,
\]
a mapping-cone resolution yields
\begin{equation}\label{eq:intro-cone}
\Ext_R^q(R/I^n,M)\cong
\Ext_S^q(S/J,M)\oplus\Ext_S^{q-1}(S/J,M).
\end{equation}
The transition maps are $(a,b)\mapsto(a,0)$, hence
\[
H_I^q(M)\cong H_{I,\A}^q(M)\cong\Ext_S^q(S/J,M).
\]
The resolution is a standard homological construction. What we make explicit is its compatibility with the entire ideal-power system, the resulting canonical projections, and the cyclic modules and comparison failures obtained from it.

Eisenbud, Musta\c{t}\u{a}, and Stillman ask when local cohomology of a polynomial ring is an increasing union of canonically embedded Ext modules and when the first canonical Ext map is injective \cite{EisenbudMustataStillman2000}. Singh and Walther obtain injectivity results using pure morphisms \cite{SinghWalther2007}; Ma, Schwede, and Shimomoto obtain characteristic-zero results related to Du Bois singularities \cite{MaSchwedeShimomoto2017}. Those are results for the original ring-valued questions. Our applications instead compute module-valued maps over non-Noetherian rings: by~\eqref{eq:intro-cone}, their kernels are the shifted Ext summands. The categorical obstruction above excludes ring-valued specialization for non-idempotent ideals in polynomial rings. Proposition~\ref{thm:union} records the direct-limit criterion used in these computations.

\subsection{Contributions and organization}
The results are organized around the following contributions.
\begin{enumerate}[label=(\roman*)]
\item A cokernel restriction on exact full abelian subcategories of reduced modules, yielding an idempotent decomposition, ambient vanishing, and finite generation for finite Noetherian coefficients; under submodule and quotient closure, vanishing extends to arbitrary coefficients.
\item An explicit polynomial-lifting comparison and a determination of the canonical Ext transition maps, giving nonzero cyclic ambient cohomology for valuation-domain and product-of-fields coefficient ideals.
\item Examples with cyclic coefficients where all ambient cohomology modules are finite and infinitely many positive-degree groups are nonzero, while relative positive-degree cohomology vanishes.
\item Applications describing kernels of canonical Ext maps and a precise obstruction to using this framework for the original polynomial-ring questions.
\end{enumerate}
The comparison machinery is recalled in Section~2; the structural results and acyclic Ext model occupy Sections~3 and~4. Section~\ref{sec:lift} contains the polynomial computations, Section~\ref{sec:applications} the Ext-map applications, and the final section the comparison failures. The contribution lies in these restrictions and explicit computations, with the standard derived-functor machinery serving as their foundation.

\section{Relative cohomology and comparison}
\subsection{The relative construction}
Throughout, $R$ is a commutative ring with identity and $I$ is an ideal. Unless otherwise stated, $\A$ is a full abelian subcategory of $R\Mod$ with enough injectives, its inclusion $j$ is exact, and all its objects are $I$-reduced. We use ambient local cohomology to mean $R^q\Gamma_I$ also when $R$ is non-Noetherian. In that generality it need not agree with geometric cohomology with supports.

If $I^n m=0$ and $n\geq2$, reducedness applied to $I^{n-2}m$ gives $I^{n-1}m=0$. Induction therefore gives the natural identification
\begin{equation}\label{eq:reduced}
\Gamma_I(M)=(0:_M I)\cong\Hom_R(R/I,M)\qquad(M\in\A).
\end{equation}
Define the $R\Mod$-valued relative functors by
\[
H_{I,\A}^q(M):=R^q_{\A}(\Gamma_I\circ j)(M)
=H^q(\Gamma_I(jE^\bullet)),
\]
where $E^\bullet$ is an injective resolution in $\A$. No closure under $\Gamma_I$ is needed for this definition. If that closure holds, the same construction can be performed internally and then embedded in $R\Mod$.

Equation~\eqref{eq:reduced} gives
\begin{equation}\label{eq:annihilation}
I H_{I,\A}^q(M)=0\qquad(q\geq0).
\end{equation}
It does not by itself identify relative cohomology with ambient Ext, since the resolutions differ.

\subsection{The comparison theorem}
Exactness of $j$ and universality of relative derived functors give canonical maps
\[
\eta_M^q:H_{I,\A}^q(M)\longrightarrow H_I^q(jM).
\]
\begin{proposition}[Comparison criterion]\label{Comp:Thm}
The maps $\eta_M^q$ are isomorphisms for every
$M\in\A$ and every $q\geq 0$ if and only if every
injective object $E$ of $\A$ is ambient
$\Gamma_I$-acyclic:
\(H_I^p(jE)=0\qquad(p>0).\)
\end{proposition}

\begin{proof}
Suppose first that the acyclicity condition holds.
Choose an injective resolution $E^\bullet$ of $M$ in $\A$.
Since $j$ is exact, $jE^\bullet$ is an ambient resolution
of $jM$. Its terms are $\Gamma_I$-acyclic by assumption.
The acyclic-resolution theorem \cite[Section~2.5]{Weibel1994} gives
\[
H_I^q(jM)
\cong H^q\bigl(\Gamma_I(jE^\bullet)\bigr)
=H^q_{I,\A}(M).
\]
To identify the map, choose an ambient injective resolution $U^\bullet$ of $jM$. The augmentation extends to a chain map $jE^\bullet\to U^\bullet$ lifting $\operatorname{id}_{jM}$; the comparison theorem makes its homotopy class independent of the choice. Applying $\Gamma_I$ gives the acyclic-resolution isomorphism above. These maps are natural in $M$ and form a morphism of cohomological $\delta$-functors extending the identity in degree zero. Universality of the relative derived functors therefore identifies them with $\eta_M^q$.

Conversely, suppose all comparison maps are isomorphisms.
If $E$ is injective in $\A$, then
\(H^p_{I,\A}(E)=0\text{  }(p>0).\) It follows that \(H_I^p(jE)\cong H^p_{I,\A}(E)=0\text{  }(p>0).\)
\end{proof}

\subsection{Quotient categories and a Tor criterion}
\begin{proposition}\label{thm:main}
Let $K$ be an ideal, $S=R/K$, and $J=(I+K)/K$. If $J^2=J$, then $\A=\{M:KM=0\}\simeq S\Mod$ satisfies the standing assumptions and
\begin{equation}\label{eq:quotient}
H_{I,\A}^q(M)\cong\Ext_S^q(S/J,M).
\end{equation}
If $0\to J\to S\to S/J\to0$ is nonsplit, then $H_{I,\A}^1(J)\ne0$.
\end{proposition}
\begin{proof}
Restriction of scalars identifies $S\Mod$ with the stated full subcategory and is exact. The $I$-action factors through $J$, so $J^2=J$ gives reducedness and identifies restricted torsion with $\Hom_S(S/J,-)$. Deriving in $S\Mod$ proves~\eqref{eq:quotient}. The nonsplit sequence represents a nonzero class in $\Ext_S^1(S/J,J)$.
\end{proof}
\begin{corollary}[Comparison for the quotient construction]
\label{cor:quotient-comparison}
Use the hypotheses and notation of Proposition~\ref{thm:main}, and
write
\(S=R/K\), \(J=(I+K)/K\), \(\mathcal A_I
=\{M\in R\text{-}\mathrm{Mod}:KM=0\}.\) Suppose further that, for every $q\geq1$,
the inverse system \(\{\operatorname{Tor}_q^R(S,R/I^n)\}_{n\geq1}\)
is pro-zero: for each $n$ there is $m\geq n$ for which the transition map to stage $n$ is zero.
Then every injective object of $\mathcal A_I$
is ambient $\Gamma_I$-acyclic. Consequently,
the canonical comparison maps of
Proposition~\ref{Comp:Thm} are isomorphisms,
and there are natural isomorphisms
\[
H_I^q(M)
\cong H_{I,\mathcal A_I}^q(M)
\cong \operatorname{Ext}_S^q(S/J,M)\text{, for every  }M\in\mathcal A_I \text{ and } q\geq0.\]
\end{corollary}

\begin{proof}
An internal injective $E$ is injective over $S$. For a projective resolution $P_\bullet\to N$, adjunction gives
\[
\Hom_R(P_\bullet,E)\cong\Hom_S(S\otimes_RP_\bullet,E).
\]
Exactness of $\Hom_S(-,E)$ yields natural isomorphisms
\[
\Ext_R^q(N,E)\cong\Hom_S(\Tor_q^R(S,N),E).
\]
Apply this to $N=R/I^n$. The Ext transition maps are induced contravariantly by the Tor transition maps. By pro-zero, for each $n$ there is $m\geq n$ for which the Ext transition from $n$ to $m$ is zero. Formula~\eqref{eq:limit} therefore gives $H_I^q(E)=0$ for $q>0$. Proposition~\ref{Comp:Thm} and Proposition~\ref{thm:main} finish the proof.
\end{proof}

\section{Structural restrictions and Noetherian vanishing}\label{sec:structural}
\subsection{A cokernel restriction}
The next observations use fullness, abelianness, and exact inclusion; they do not require enough injectives.
\begin{proposition}\label{prop:coreduced}
Suppose $I^2$ is finitely generated and $\mathcal B\subseteq R\Mod$ is a full abelian subcategory with exact inclusion whose objects are $I$-reduced. Then
\[
IM=I^2M\qquad(M\in\mathcal B).
\]
In particular, $R\in\mathcal B$ implies $I^2=I$.
\end{proposition}
\begin{proof}
Choose generators $a_1,\ldots,a_s$ of $I^2$. Since $\mathcal B$ is additive and full, the map
\[
\phi:M^s\to M,\qquad (m_1,\ldots,m_s)\mapsto\sum_i a_i m_i
\]
is a morphism in $\mathcal B$. Its internal cokernel has underlying module $M/I^2M$, because the inclusion is exact. This quotient is $I$-reduced and killed by $I^2$, so it is killed by $I$. Hence $IM\subseteq I^2M$, proving equality. For $M=R$ this is $I=I^2$.
\end{proof}

\subsection{Finite coefficients over Noetherian rings}
\begin{theorem}\label{thm:finite-vanishing}
Let $R$ be Noetherian, and let $\mathcal B$ satisfy the categorical hypotheses of Proposition~\ref{prop:coreduced}. For every finitely generated $M\in\mathcal B$,
\[
I_{\mathfrak p}M_{\mathfrak p}=0\quad(\mathfrak p\supseteq I),
\qquad
H_I^q(M)=0\quad(q>0).
\]
Consequently, $H_I^q(M)$ is finitely generated for every $q\geq0$.
\end{theorem}
\begin{proof}
Proposition~\ref{prop:coreduced} gives $IM=I^2M$. At a prime $\mathfrak p\supseteq I$, the finitely generated $R_{\mathfrak p}$-module $N=I_{\mathfrak p}M_{\mathfrak p}$ satisfies $N=I_{\mathfrak p}N$. Nakayama's lemma gives $N=0$. Over a Noetherian ring, local cohomology commutes with localization and is computed by an augmented \v{C}ech complex on generators of $I$. At primes containing $I$, the positive terms of that complex vanish because $I_{\mathfrak p}M_{\mathfrak p}=0$. At other primes the supporting ideal becomes the unit ideal and local cohomology is zero. Thus $H_I^q(M)$ has zero localization at every prime for $q>0$. The degree-zero module is a submodule of the finite module $M$ and is therefore finite.
\end{proof}

\begin{corollary}[Decomposition on the effective coefficient ring]\label{cor:decomposition}
Under the hypotheses of Theorem~\ref{thm:finite-vanishing}, put $S=R/\Ann_R(M)$. There is an idempotent $e\in S$ such that the image of $I$ is $eS$, and
\[
M=\Gamma_I(M)\oplus IM,\qquad \Gamma_I(M)=(1-e)M,\quad IM=eM.
\]
\end{corollary}
\begin{proof}
By Proposition~\ref{prop:coreduced}, the finite module $N=IM$ satisfies $N=IN$. The determinant trick gives $r\in I$ with $(1-r)N=0$. Thus $(r-r^2)M=0$, and $e=r+\Ann_R(M)$ is idempotent. For $a\in I$, $(1-r)aM=0$, so the image of $a$ lies in $eS$. Conversely $e$ belongs to the image of $I$, giving equality. The idempotent decomposition of $M$ now applies. The ideal $I$ kills $(1-e)M$, whereas any element of $eM$ killed by a power of $I$ is killed by a power of $e$ and hence is zero. This identifies the torsion summand.
\end{proof}

\begin{remark}\label{rem:finite-scope}
The theorem does not require enough internal injectives or either acyclicity assumption used below. In particular, a finite-generation statement for ambient cohomology of finite objects in this Noetherian framework gives vanishing in positive degrees. The nonzero finite-generation examples in Section~\ref{sec:lift} instead use non-Noetherian rings and non-finitely-generated coefficient ideals.
\end{remark}

\subsection{Vanishing for arbitrary objects under closure}
The closure conditions in the following theorem force reducedness to hold on the cyclic subquotients needed for a local argument. This proves ambient vanishing without assuming comparison in advance.
\begin{theorem}\label{thm:vanish}
Let $R$ be Noetherian. Suppose $\mathcal A$ is a full abelian subcategory of $R\Mod$, closed under $R$-submodules and quotients, whose objects are $I$-reduced. Then $\Gamma_I|_{\mathcal A}$ is exact and
\[
H_I^q(M)=0\qquad(M\in\mathcal A,\ q>0).
\]
If $\mathcal A$ has enough injectives, its relative cohomology also vanishes in positive degrees and the canonical comparison maps are isomorphisms. For finitely generated $M\in\mathcal A$, $H_I^q(M)$ is finitely generated in every degree.
\end{theorem}
\begin{proof}
For $m\in M$, put $C=Rm$. Both $C$ and $C/I^2C$ belong to $\mathcal A$. Reducedness of $C/I^2C$, which is killed by $I^2$, gives $IC=I^2C$. At a prime $\mathfrak p\supseteq I$, the finitely generated module $I_{\mathfrak p}C_{\mathfrak p}$ therefore satisfies
\[
I_{\mathfrak p}C_{\mathfrak p}=I_{\mathfrak p}(I_{\mathfrak p}C_{\mathfrak p}).
\]
Nakayama's lemma gives $I_{\mathfrak p}C_{\mathfrak p}=0$. Since $m$ was arbitrary,
\begin{equation}\label{eq:localann}
I_{\mathfrak p}M_{\mathfrak p}=0\qquad(\mathfrak p\supseteq I).
\end{equation}
Local cohomology commutes with localization over a Noetherian ring. At primes containing $I$, the positive terms of a \v{C}ech complex vanish by~\eqref{eq:localann}; at the other primes the supporting ideal becomes the unit ideal and the complex is acyclic. Thus every localization of $H_I^q(M)$ is zero for $q>0$.

Fullness and closure under submodules and quotients ensure that ambient kernels and cokernels of morphisms in $\mathcal A$ belong to $\mathcal A$ and satisfy their internal universal properties. Thus the inclusion is exact. For a short exact sequence in $\mathcal A$, the ambient local cohomology sequence and the vanishing of $H_I^1$ of its first term show that restricted torsion preserves surjections. It is therefore exact. Enough internal injectives then gives relative vanishing. All objects, in particular internal injectives, are ambient torsion-acyclic, so Proposition~\ref{Comp:Thm} proves comparison. Finally, $H_I^0(M)\subseteq M$ is finitely generated when $M$ is; higher groups are zero.
\end{proof}

\begin{example}
Take $S=k\times k$, $J=k\times0=eS$, $R=S[t]$, $I=JR+(t)$, and $\A=\{M:tM=0\}$. This subcategory satisfies the theorem. The ideal $I$ is non-idempotent, while
\[
H_I^q(M)\cong H_{I,\A}^q(M)\cong
\begin{cases}(1-e)M,&q=0,\\0,&q>0.\end{cases}
\]
Thus the supporting ideal need not be idempotent for its torsion functor to become an exact nonzero direct-summand functor on the subcategory.
\end{example}

\begin{remark}
More generally, in the quotient construction over a Noetherian ring, $J$ is a finitely generated idempotent ideal and hence $J=eS$. Then $S/J$ is projective and relative cohomology vanishes in positive degrees. Nonsplit idempotent-ideal constructions necessarily lie outside the Noetherian setting.
\end{remark}
\section{Acyclic resolutions and an ambient Ext model}
The next statements record standard acyclic-resolution consequences needed for the applications. They do not assert positive-degree nonvanishing for finite coefficients over a Noetherian ring; Theorem~\ref{thm:finite-vanishing} rules out such nonvanishing under the standing categorical assumptions.

\begin{proposition}[An ambient Ext model]\label{MainThm}
Under the standing assumptions, suppose each injective object $E$ of $\A$ satisfies
\begin{equation}\label{eq:Extacyclic}
\Ext_R^p(R/I,E)=0\qquad(p>0).
\end{equation}
Then $H_{I,\A}^q(M)\cong\Ext_R^q(R/I,M)$ naturally. If, in addition, these injectives are ambient $\Gamma_I$-acyclic, then
\begin{equation}\label{eq:triple}
\Ext_R^q(R/I,M)\xrightarrow{\ \sim\ }H_I^q(M)\cong H_{I,\A}^q(M),
\end{equation}
where the first arrow is the canonical map. If $R$ is Noetherian and $M$ is finitely generated, all these modules are finitely generated.
\end{proposition}
\begin{proof}
An internal injective resolution $E^\bullet$ is an exact ambient resolution. Under~\eqref{eq:Extacyclic} it is acyclic for $\Hom_R(R/I,-)$, so the acyclic-resolution theorem and~\eqref{eq:reduced} identify its Hom cohomology with both ambient Ext and relative cohomology. Under ambient torsion-acyclicity the same resolution computes $H_I^q(M)$. The canonical natural transformation $\Hom_R(R/I,-)\to\Gamma_I$ is a termwise isomorphism on $E^\bullet$; hence the isomorphism in~\eqref{eq:triple} is the canonical one. Finally, over a Noetherian ring, Ext between finitely generated modules is finitely generated.
\end{proof}

\begin{proposition}\label{Inj n}
Suppose $R$ is Noetherian and~\eqref{eq:Extacyclic} holds. If each $I^n/I^{n+1}$, $n\geq1$, is projective over $R/I$, every internal injective is ambient $\Gamma_I$-acyclic. Thus all conclusions of Proposition~\ref{MainThm} apply.
\end{proposition}
\begin{proof}
Each $I^n/I^{n+1}$ is a summand of a finite free $R/I$-module. Its positive Ext groups into an internal injective $E$ vanish by~\eqref{eq:Extacyclic}. The exact sequences
\[
0\to I^n/I^{n+1}\to R/I^{n+1}\to R/I^n\to0
\]
then give $\Ext_R^p(R/I^n,E)=0$ for all $n,p\geq1$ by induction. Taking the direct limit yields $H_I^p(E)=0$.
\end{proof}

These hypotheses are sufficient conditions, not a consequence of reducedness alone. In particular, the Ext-acyclicity condition must be checked independently in an application.

\begin{example}[Acyclicity hypotheses are distinct]
Let $R=k[t]$, $I=(t)$, and $\A=\{M:tM=0\}\simeq k\Mod$. Restricted torsion is exact, and every object is ambient torsion-acyclic by the two-term \v{C}ech complex. Nevertheless, the internal injective $k$ satisfies $\Ext_R^1(R/I,k)\cong k\ne0$. Thus comparison need not provide the ambient Ext model in Proposition~\ref{MainThm}.

For an example satisfying both hypotheses, let $R=A\times B\times C$ and $I=A\times0\times L$, where $L$ is an ideal of $C$. Take the full subcategory of modules $(X,Y,0)$. Projection to the first two ring factors is a flat quotient, so its restriction-of-scalars functor preserves injectives: its left adjoint is exact. Every object of this subcategory is $I$-reduced, and the internal injectives are ambient injectives. Both acyclicity hypotheses follow. Here ambient and relative torsion are $(X,Y,0)\mapsto(0,Y,0)$, and the positive cohomology groups vanish. Choosing $L^2\ne L$ makes $I$ non-idempotent.
\end{example}

\section{Polynomial lifting and non-Noetherian examples}\label{sec:lift}
The quotient construction also gives nonzero finitely generated ambient cohomology without Noetherianity. The next result supplies comparison for all coefficient modules in a polynomial lift.
\begin{theorem}\label{thm:lift}
Let $J\subsetneq S$ be idempotent, set $T:=S/J$, $R:=S[t]$, $I:=JR+(t)$, and $\A:=\{M:tM=0\}$. Then $I^2\ne I$. For every $M\in\A$, $n\geq1$, and $q\geq0$, there are natural identifications
\[
\Ext_R^q(R/I^n,M)\cong\Ext_S^q(T,M)\oplus\Ext_S^{q-1}(T,M),
\]
where negative Ext groups are zero. The maps induced by $R/I^{n+1}\twoheadrightarrow R/I^n$ are $(a,b)\mapsto(a,0)$. Consequently, there are natural isomorphisms
\[
H_I^q(M)\cong H_{I,\A}^q(M)\cong\Ext_S^q(T,M).
\]
Under these identifications, the canonical map from each Ext stage to $H_I^q(M)$ is projection onto the first summand.
\end{theorem}
\begin{proof}
Idempotence gives $I^n=JR+(t^n)$, hence $R/I^n=T[t]/(t^n)$. The images of $I$ and $I^2$ in $T[t]$ are $(t)$ and $(t^2)$, and differ since $T\ne0$.

Choose a projective resolution $P_\bullet\to T$ over $S$, concentrated in nonnegative homological degrees. Set $Q_\bullet=R\otimes_SP_\bullet$. Flatness of $R$ over $S$ makes $Q_\bullet$ a projective resolution of $T[t]$. Let $C^{(n)}_\bullet$ be the cone of $t^n:Q_\bullet\to Q_\bullet$, with
\[
C^{(n)}_i=Q_i\oplus Q_{i-1},\qquad
 d^{(n)}(a,b)=(d_Qa+t^nb,-d_Qb),
\]
where $Q_{-1}=0$. The cone homology sequence and injectivity of $t^n$ on $T[t]$ show that $C^{(n)}_\bullet$ resolves $T[t]/(t^n)$.

For $tM=0$, precomposition by $t^n$ is zero on $\Hom_R(Q_\bullet,M)$. Adjunction therefore splits the Hom complex into the complex $\Hom_S(P_\bullet,M)$ and its degree-one shift, up to the conventional sign in the shifted differential. Thus
\begin{equation}\label{eq:cone}
\Ext_R^q(R/I^n,M)\cong\Ext_S^q(T,M)\oplus\Ext_S^{q-1}(T,M),
\end{equation}
where negative Ext groups are zero.

The chain map
\[
f_n:C^{(n+1)}_\bullet\to C^{(n)}_\bullet,
\qquad f_n(a,b)=(a,tb),
\]
lifts the canonical quotient map between the modules being resolved. The chain-map identity is $d^{(n)}f_n=f_nd^{(n+1)}$. Applying the Hom functor, its map on the two summands is $(a,b)\mapsto(a,0)$. Taking the direct limit gives $H_I^q(M)\cong\Ext_S^q(T,M)$ naturally in $M$.

If $E$ is $S$-injective, these groups vanish for $q>0$. Proposition~\ref{thm:main} and Proposition~\ref{Comp:Thm} therefore identify this ambient cohomology with the relative construction. The same cone calculation identifies the canonical maps from every Ext stage with the projections onto the first summand.
\end{proof}

\begin{lemma}\label{lem:endomorphism}
For an ideal $J\subseteq S$,
\begin{equation}\label{eq:first}
\Ext_S^1(S/J,J)\cong\operatorname{End}_S(J)/\rho(J),
\qquad \rho(a)(b)=ab.
\end{equation}
If multiplication identifies $S$ with $\operatorname{End}_S(J)$ and $J^2=J$, the lift of Theorem~\ref{thm:lift} satisfies $H_I^1(J)\cong S/J\cong R/I$.
\end{lemma}
\begin{proof}
Apply $\Hom_S(-,J)$ to $0\to J\to S\to S/J\to0$. Since $S$ is projective, the resulting sequence ends in
\[
\Hom_S(S,J)\longrightarrow\operatorname{End}_S(J)
\longrightarrow\Ext_S^1(S/J,J)\longrightarrow0.
\]
Under $\Hom_S(S,J)\cong J$, the first map is $\rho$. The remaining assertion follows from Theorem~\ref{thm:lift}.
\end{proof}

The next examples realize this cyclic module. Both the coefficient ideals and the supporting ideals in these examples are infinitely generated. Thus these are examples of finite generation for derived ideal-power torsion outside the Noetherian finite-coefficient setting.
\begin{example}[A nondiscrete valuation domain]\label{ex:valuation}
Let $k$ be a field and consider the nested union \(S:=\displaystyle\bigcup_{n\geq1}k[[u^{1/n!}]],\) where the inclusions are determined by $u^{1/n!}= (u^{1/(n+1)!})^{n+1}$. Its fraction field is \(F:=\displaystyle\bigcup_{n\geq1}k((u^{1/n!})).\) The usual valuations on these Laurent series fields are compatible and define a valuation $v:F^\times\to\mathbb Q$, normalized by $v(u)=1$, and
\(S=\{0\}\cup\{a\in F^\times:v(a)\geq0\}.\) Thus $S$ is a valuation domain with maximal ideal \(\mathfrak m=\{0\}\cup\{a\in F^\times:v(a)>0\}\text{, } S/\mathfrak m\cong k.\)
For $0\neq a\in\mathfrak m$, put $q=v(a)>0$. Then \(a=u^{q/2}(u^{-q/2}a),\) and both factors have valuations $q/2$, so $a\in \mathfrak m^2$. Thus $\mathfrak m\subseteq\mathfrak m^2$; the reverse inclusion is automatic, hence $\mathfrak m^2=\mathfrak m$. Since $S$ is a domain and $0\neq\mathfrak m\subsetneq S$, the ideal $\mathfrak m$ cannot be generated by an idempotent. The characterization of split ideals by idempotent generators therefore shows that its canonical sequence \(  0\longrightarrow \mathfrak m \xrightarrow{\iota}S \longrightarrow S/\mathfrak m \longrightarrow 0\) is nonsplit. Set $R=S[x]$, $I=\mathfrak mS[x]+(x)$ and
$\mathcal A_I=\{M:xM=0\}$. Theorem~\ref{thm:lift} applies. We now compute the first ambient and relative cohomology of $\mathfrak m$,
regarded as an $R$-module via $R/(x)\cong S$. Every $S$-endomorphism of the nonzero rank-one torsion-free module
$\mathfrak m$ is multiplication by an element of $F$. Therefore
\(\operatorname{End}_S(\mathfrak m)\cong(\mathfrak m:\mathfrak m)=\{c\in F:c\mathfrak m\subseteq\mathfrak m\}.\) Clearly $S\subseteq(\mathfrak m:\mathfrak m)$. Conversely, if $v(c)<0$, choose $r\in\mathbb Q$ with $0<r<-v(c)$. Then $u^r\in\mathfrak m$ but $cu^r\notin\mathfrak m$.
Thus $(\mathfrak m:\mathfrak m)=S$, and multiplication induces an isomorphism $S\cong\operatorname{End}_S(\mathfrak m)$. By Theorem~\ref{thm:lift} and Equation~\eqref{eq:first}, \[H^1_I(\mathfrak m)\cong\operatorname{Ext}_S^1(S/\mathfrak m,\mathfrak m)\cong S/\mathfrak m\cong k\neq0.\]
In particular, this cohomology module is cyclic over $R$. The ideal $\mathfrak m$ is not finitely generated: in a valuation domain a finitely generated ideal is principal, whereas the positive rational values have no least element. The image of $I$ under evaluation at $x=0$ is $\mathfrak m$, so $I$ is not finitely generated either.
\end{example}

\begin{example}[A product of fields]\label{ex:product}
Let $S:=\prod_{n\geq1}k$, $J:=\bigoplus_{n\geq1}k$, and $T:=S/J$. \(J\) is regarded as an \(R:=S[t]\)-module through evaluation at \(t=0\). If $e_n$ denotes the coordinate idempotent, then $J=\bigoplus_n Se_n$ is projective and $J^2=J$. Every endomorphism preserves the summands, so multiplication gives $\operatorname{End}_S(J)\cong S$. The ideal is not generated by an idempotent: any such generator in $J$ has finite support, whereas $J$ contains every $e_n$.
For the lift $R=S[t]$, $I=JR+(t)$, Theorem~\ref{thm:lift} and~\eqref{eq:first} yield
\[
H_I^q(J)\cong H_{I,\A}^q(J)\cong
\begin{cases}T,&q=1,\\0,&q\ne1.\end{cases}
\]
Here $\Hom_S(T,J)=0$ and $\operatorname{pd}_S T\leq1$. In particular, $H_I^1(J)\cong R/I$ is nonzero and cyclic. Formula~\eqref{eq:cone} also gives $\Ext_R^2(R/I^n,J)\cong T$ at every stage, but all its transition maps are zero. Thus these degree-two Ext classes contribute nothing to ambient local cohomology. The ideal $J$ is not finitely generated, since finitely many elements have support in one fixed finite set. Evaluation at $t=0$ maps $I$ onto $J$, so $I$ is also not finitely generated.
\end{example}
\section{Applications}\label{sec:applications}
\subsection{The questions of Eisenbud, Musta\c{t}\u{a}, and Stillman}
For a polynomial ring $P=k[x_1,\ldots,x_r]$, Eisenbud, Musta\c{t}\u{a}, and Stillman ask which ideals admit a descending sequence $I\supseteq L_1\supseteq L_2\supseteq\cdots$ such that the canonical Ext maps are injective and their images exhaust $H_I^q(P)$, and which ideals have injective canonical map
\[
\Ext_P^q(P/I,P)\longrightarrow H_I^q(P);
\]
see \cite[Questions~6.1 and~6.2]{EisenbudMustataStillman2000}. A standard sufficient setting for the exhaustion is to choose $(L_n)$ cofinal with the powers of $I$. We use such cofinal systems below; cofinality is an additional condition in this formulation, rather than a replacement for the original exhaustion question.

Their monomial results and the subsequent work \cite{MaSchwedeShimomoto2017,SinghWalther2007} concern the original ring-valued problems. Our conclusions concern coefficient modules $M$ over general rings. Proposition~\ref{prop:coreduced} shows that the present categorical framework cannot contain $P$ when $I$ is non-idempotent. Accordingly, the following statements are module-valued criteria and examples, not new solutions of those original polynomial-ring questions.

\subsection{A criterion for an increasing union of Ext-submodules}
\begin{proposition}[A direct-limit criterion]\label{thm:union}
Let $R$ be any commutative ring, $I$ an ideal, $M$ an $R$-module, and $q\geq1$. Write
\[
\alpha_n^q:\Ext_R^q(R/I^n,M)\to H_I^q(M).
\]
The following are equivalent:
\begin{enumerate}[label=(\roman*)]
\item every $\alpha_n^q$ is injective;
\item every transition $\Ext_R^q(R/I^n,M)\to\Ext_R^q(R/I^{n+1},M)$ is injective;
\item every boundary map
$\Ext_R^{q-1}(I^n/I^{n+1},M)\to\Ext_R^q(R/I^n,M)$ is zero.
\end{enumerate}
Under these conditions $H_I^q(M)$ is the increasing union of these canonically embedded Ext modules. In particular, vanishing of all the displayed source groups is sufficient. In degree zero the canonical maps are always injective.
\end{proposition}
\begin{proof}
The direct-limit formula identifies $\alpha_n^q$ with the structure map of the Ext system. Injectivity of every structure map implies injectivity of every transition. Conversely, if every transition is injective, a class that becomes zero in the limit is zero at its starting stage. The long exact Ext sequence associated to
\[
0\to I^n/I^{n+1}\to R/I^{n+1}\to R/I^n\to0
\]
identifies the kernel of the transition map with the image of the stated boundary map. Every direct-limit class comes from a stage, giving the union assertion. In degree zero the maps are the inclusions of annihilator submodules into $\Gamma_I(M)$.
\end{proof}

The same argument works for any descending sequence $L_n$ cofinal with the powers of $I$, replacing $I^n/I^{n+1}$ by $L_n/L_{n+1}$. Proposition~\ref{thm:union} is a general direct-limit criterion; relative cohomology enters by providing a model for the target and, under further acyclicity, identifying its canonical maps.

\subsection{First-stage maps from acyclic resolutions}
Under both acyclicity hypotheses of Proposition~\ref{MainThm}, the canonical map $\Ext_R^q(R/I,M)\to H_I^q(M)$ is an isomorphism for $M\in\A$. This does not imply that all later Ext stages embed. For the power sequence, Proposition~\ref{thm:union} specifies the additional requirement; for example, $\Ext_R^{q-1}(I^n/I^{n+1},M)=0$ for all $n\geq1$ is sufficient in degree $q\geq1$.

\subsection{Explicit answers and obstructions in polynomial lifts}
Theorem~\ref{thm:lift} makes the injectivity issue completely explicit. For its notation and $tM=0$, the canonical maps in positive degrees are, under~\eqref{eq:cone}, the projections
\[
\Ext_S^q(T,M)\oplus\Ext_S^{q-1}(T,M)
\longrightarrow\Ext_S^q(T,M).
\]
Indeed, the direct system is $(a,b)\mapsto(a,0)$ and the proof of comparison identifies its limit with the displayed target. Consequently, all the power-stage maps in degree $q$ are injective if and only if
\[
\Ext_S^{q-1}(T,M)=0.
\]
In the product-of-fields example with $M=J$, the degree-one maps are isomorphisms onto the nonzero cyclic module $T$, because $\Hom_S(T,J)=0$. Thus
\[
H_I^1(J)=\bigcup_{n\geq1}\operatorname{im}\bigl(\Ext_R^1(R/I^n,J)\to H_I^1(J)\bigr),
\]
and each map is an isomorphism. In degree two each source is $T\ne0$ while the target is zero, so injectivity fails at every stage. These are concrete affirmative and negative instances of the module-valued questions for a non-idempotent ideal.

\subsection{The obstruction to ring-valued specialization}
Proposition~\ref{prop:coreduced} implies that, whenever $I^2$ is finitely generated, a full abelian subcategory of $I$-reduced modules with exact inclusion can contain $R$ only if $I^2=I$. No closure under arbitrary quotient modules is needed for this conclusion. In a polynomial ring over a field every ideal is finitely generated, so the obstruction applies to every non-idempotent supporting ideal. The module-valued examples above therefore cannot be specialized to the original ring-valued questions by simply taking $M=R$.

\section{Why comparison cannot be omitted}
The next examples show that even a coefficient module annihilated by the supporting ideal can have nonzero positive-degree ambient derived torsion over a non-Noetherian ring. Its relative cohomology in a semisimple quotient category nevertheless vanishes.
\subsection{An idempotent supporting ideal}\label{subsec:counterexample}

Let
\[
A=
k[u_1,u_2,\ldots]\Big/
\bigl(u_1^2,\ u_n-u_{n+1}^2\ (n\geq 1)\bigr),
\]
and put
\(\mathfrak m=(u_1,u_2,\ldots).\)
The relations give \(\mathfrak m^2=\mathfrak m, \text{ and } A/\mathfrak m\cong k.\)
Take $I=\mathfrak m$ and \(\A=\{M:\mathfrak mM=0\}\simeq k\Mod.\)
This is a full abelian subcategory with exact inclusion
and enough injectives. Its objects are $I$-reduced. Since restricted torsion is the identity on this category,
\[H^q_{I,\A}(k)=
\begin{cases}
k,&q=0,\\
0,&q>0.
\end{cases}
\]

In the ambient category, $I^n=I$ for every $n\geq 1$,
so
\[
\Gamma_I(-)=\Hom_A(k,-),
\qquad
H_I^q(k)=\Ext_A^q(k,k).
\]

\begin{proposition}\label{prop:even-ext}
For the ring $A$ above,
\[
\Ext_A^q(k,k)\cong
\begin{cases}
k,&q\text{ even},\\
0,&q\text{ odd}.
\end{cases}
\]
\end{proposition}

\begin{proof}
Write
\[
A_n=k[u_n]/(u_n^{2^n}),
\qquad
A=\varinjlim_n A_n,
\]
where the transition map sends $u_n$ to $u_{n+1}^2$. These maps are injective, and $A_{n+1}$ is free over $A_n$ with basis $1,u_{n+1}$. The residue field $k$ has an $A_n$-free resolution
\[
\cdots
\xrightarrow{u_n}A_n
\xrightarrow{u_n^{2^n-1}}A_n
\xrightarrow{u_n}A_n
\longrightarrow k\longrightarrow 0.
\]
Tensoring with $k$ makes all differentials zero, giving
\(\Tor_q^{A_n}(k,k)\cong k
\qquad(q\geq 0).\)

A semilinear chain map corresponding to $A_n\to A_{n+1}$
multiplies by $1$ in even degrees and by $u_{n+1}$
in odd degrees. To verify compatibility, write
$v=u_{n+1}$ and $a=2^n$; then $u_n$ maps to $v^2$,
and the identities needed for the chain map are
\[
v\cdot v=v^2,
\qquad
v^{2a-1}\cdot 1=v\cdot v^{2a-2}.
\]
After tensoring with $k$, the induced Tor maps are
identities in even degrees and zero in odd degrees.

For an augmented $k$-algebra $B$ with augmentation ideal $\overline B$, the reduced bar complex computing $\Tor^B(k,k)$ has degree-$q$ term $\overline B^{\otimes_k q}$ and differential induced by multiplication. See \cite[Chapter~8]{Weibel1994}. This construction is functorial in augmented algebra maps. Tensor products and filtered colimits commute, and filtered colimits of vector spaces are exact. Thus the bar complexes identify
\[
\Tor_q^A(k,k)\cong\varinjlim_n\Tor_q^{A_n}(k,k).
\]
To justify the identification of maps, regard the $A_{n+1}$-resolution as an exact complex of $A_n$-modules. The displayed semilinear map is an $A_n$-linear lift of the identity on $k$. Projectivity of the source resolution implies that any two such lifts are chain homotopic. Comparing with the functorial bar resolutions shows that its induced map is exactly the canonical change-of-rings map on Tor. Hence
\[
\Tor_q^A(k,k)\cong
\begin{cases}k,&q\text{ even},\\0,&q\text{ odd}.\end{cases}
\]
The free bar resolution also identifies the cochain complex computing $\Ext_A^*(k,k)$ with the $k$-linear dual of the bar chain complex computing $\Tor_*^A(k,k)$. Vector-space duality is exact, even in infinite dimension, so
\[
\Ext_A^q(k,k)\cong\Hom_k(\Tor_q^A(k,k),k).
\]
This proves the assertion.
\end{proof}

It follows that
\[H^2_{I,\A}(k)=0,\qquad H_I^2(k)\cong k\neq 0.\]
More generally, the two theories differ in every
positive even degree.

\subsection{A non-idempotent supporting ideal}

\begin{theorem}[Cyclic coefficients with unbounded ambient cohomology]\label{thm:unbounded}
Let $A$ and $\mathfrak m$ be as in Proposition~\ref{prop:even-ext}, and put
\[
R=A[t],\qquad I=\mathfrak mR+(t),\qquad
\A=\{M:IM=0\}\simeq k\Mod.
\]
The ideal $I$ is non-idempotent and infinitely generated. The category $\A$ is a full abelian subcategory of $I$-reduced modules with exact inclusion and enough injectives. The cyclic coefficient module $k=R/I$ satisfies
\[
H_{I,\A}^q(k)=0\quad(q>0),\qquad
H_I^q(k)\cong
\begin{cases}k,&q\text{ even},\\0,&q\text{ odd}.\end{cases}
\]
Thus, every ambient cohomology module is finitely generated, and comparison fails in every positive even degree.
\end{theorem}
\begin{proof}
The image of $I$ in $R/\mathfrak mR\cong k[t]$ is $(t)$, so $I^2\ne I$. Restriction of scalars identifies $\A$ with $k\Mod$, giving its asserted categorical properties. Restricted torsion is the identity, so its positive derived functors vanish.

For the ambient computation,
\(I^n=\mathfrak mR+(t^n)\text{, } R/I^n\cong k[t]/(t^n).\)
Extend an $A$-free resolution of $k$ to $A[t]$ and
take the mapping cone of multiplication by $t^n$.
Since $tk=0$, this gives
\(\Ext_R^q(R/I^n,k)\cong\Ext_A^q(k,k)\oplus\Ext_A^{q-1}(k,k).\)
The transition maps are identities on the first summand
and zero on the second. Consequently,
\[
H_I^q(k)
\cong
\varinjlim_n\Ext_R^q(R/I^n,k)
\cong
\Ext_A^q(k,k).
\]

We therefore obtain
\[
H_I^q(k)\cong
\begin{cases}
k,&q\text{ even},\\
0,&q\text{ odd}.
\end{cases}
\]
The comparison is:
\[
\begin{array}{c|c|c}
q&H^q_{I,\A}(k)&H_I^q(k)\\ \hline
0&k&k\\
\text{positive odd}&0&0\\
\text{positive even}&0&k
\end{array}
\]

The ideal $\mathfrak m$ is infinitely generated: a finite set of its elements lies in some $A_N$, so the ideal it generates is contained in $u_NA$, whereas $u_{N+1}\notin u_NA$. To see the latter, use injectivity of the transition maps and compare monomial exponents in a sufficiently large truncated polynomial ring. The polynomial-lift ideal is infinitely generated as well, since evaluation at $t=0$ maps it onto $\mathfrak m$.
\end{proof}

Note that in the polynomial construction of Section 7, there are two relevant categories:
\[
\mathcal B=\{M:tM=0\}\simeq A\text{-}\mathrm{Mod},
\text{ and }
\mathcal C=\{M:IM=0\}\simeq k\text{-}\mathrm{Mod}.
\]
Theorem 5.1 gives comparison with relative cohomology in \(\mathcal B\). Theorem 7.2 demonstrates failure of comparison in the smaller category \(\mathcal C\).

\begin{remark}[Why coincidence fails]
In both constructions, the module $k$ is injective
inside $\A\simeq k\Mod$, but \(H_I^2(k)\cong k\neq 0.\)
Thus internal injectives are not ambient
$\Gamma_I$-acyclic. This is precisely the failure
identified by the comparison theorem.
\end{remark}

\begin{remark}[Derived torsion versus geometric support]
In Subsection~\ref{subsec:counterexample}, every element of $\mathfrak m$ is
nilpotent, so \(V(\mathfrak m)=\Spec A.\)
Since \(V(\mathfrak m)=\operatorname{Spec}A\), geometric local cohomology of the associated quasi-coherent sheaf is ordinary sheaf cohomology on the affine scheme \(\operatorname{Spec}A\), and therefore vanishes in positive degrees. Nevertheless,
the derived $\mathfrak m$-power torsion functors above
are nonzero in positive even degrees. This illustrates why the distinction between derived
power torsion and geometric local cohomology must be
retained for these non-Noetherian examples.
\end{remark}

\section*{Acknowledgements}
The authors used ChatGPT (OpenAI) to assist with the revision of the manuscript. However, all research questions were formulated by the
authors, and they take full responsibility for
the mathematical results, proofs, references, and final text.
\section*{Competing interests}
The authors declare no competing interests.

\end{document}